\documentclass[11pt,a4paper]{amsart}

\usepackage[T1]{fontenc}
\usepackage[utf8]{inputenc}
\usepackage[english]{babel}
\usepackage{amsmath,amssymb,amsthm}
\usepackage{mathtools}
\usepackage{tikz}
\usepackage{geometry}
\usepackage[hidelinks]{hyperref}
\usetikzlibrary{arrows.meta,calc,fit,positioning,decorations.pathreplacing}
\theoremstyle{definition}
\newtheorem{definition}{Definition}[section]

\theoremstyle{plain}
\newtheorem{lemma}[definition]{Lemma}
\newtheorem{proposition}[definition]{Proposition}
\newtheorem{theorem}[definition]{Theorem}

\theoremstyle{remark}
\newtheorem{remark}[definition]{Remark}

\newcommand{\len}[1]{\ell(#1)}
\newcommand{\dist}[3]{d_{#1}(#2,#3)}

\title{Structural properties of semi-cutting sets}

\author{Bart{\l}omiej Bosek}
\address{Institute of Theoretical Computer Science, Faculty of Mathematics and Computer Science, Jagiellonian University, Krak\'ow, Poland}
\email{bartlomiej.bosek@uj.edu.pl}
\thanks{Bart{\l}omiej Bosek is supported by the National Science Centre, Poland, under grant no. 2023/51/B/ST6/02833.}

\author{Oriol Serra}
\address{Department of Mathematics, Universitat Polit\`ecnica de Catalunya, Barcelona, Spain and Centre de Recerca Matem\`atica, Barcelona, Spain}
\email{oriol.serra@upc.edu}
\thanks{Oriol Serra and Llu\'is Vena  were supported by the Spanish Agencia Estatal de Investigación under projects PID2020-113082GB-I00, PID2023-147202NB-I00, RED2022-134947-T and the Severo Ochoa and María de Maeztu Program for Centers and Units of Excellence in R\&D (CEX2020-001084-M)}

\author{Llu{\'{\i}}s Vena}
\address{Department of Mathematics, Universitat Polit\`ecnica de Catalunya, Barcelona, Spain}
\email{lluis.vena@upc.edu}

\date{\today}


\begin{document}

\begin{abstract}
Let $G=(V,E)$ be a finite simple undirected graph, let $T\subseteq V$, and let $F\subseteq E$. We call $F$ \emph{semi-cutting} for $T$ if $\dist{G}{u}{w}<\dist{G\setminus F}{u}{w}$ for every two distinct terminals $u,w\in T$. We prove the sharp bound $|F|\ge |T|-1$. The argument is structural: semi-cutting sets are first identified with edge transversals of all shortest terminal-terminal paths, and the bound then follows from a contraction step, a terminal-deletion lemma, and a lexicographic induction on $(|T|,|V|)$.
\end{abstract}

\maketitle

\section{Introduction}

All graphs in this paper are finite, simple, and undirected. In particular, loops and parallel edges are excluded, and every contraction is understood in the simple-graph sense: loops are deleted and parallel edges are identified. For a graph $G=(V,E)$, a terminal set $T\subseteq V$, and an edge set $F\subseteq E$, we say that $F$ is \emph{semi-cutting} for $T$ if $\dist{G}{u}{w}<\dist{G\setminus F}{u}{w}$ for every two distinct terminals $u,w\in T$, where distances take values in $\mathbb{N}_0\cup\{\infty\}$. No global connectivity assumption is imposed on $G$. Indeed, if some terminal pair already had infinite distance in $G$, then the defining inequality would read $\infty<\infty$ and could not hold.

The condition says that every shortest path whose endpoints lie in $T$ must lose at least one edge after the deletion of $F$. A lower bound of order $|T|$ is therefore natural. It is also nontrivial: one edge may lie on shortest paths for many terminal pairs, so a direct count over pairs gives no useful control. Our main result is that every semi-cutting set satisfies $|F|\ge |T|-1$.

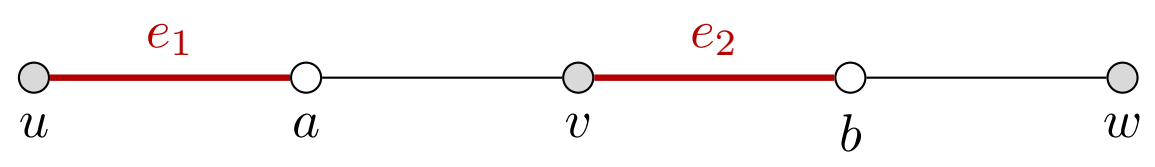
\begin{figure}[htbp]
\centering
\begin{tikzpicture}[scale=1, every node/.style={font=\small}]
  \node[circle,draw,inner sep=2pt,fill=gray!30,label=below:$u$] (u) at (0,0) {};
  \node[circle,draw,inner sep=2pt,fill=white,label=below:$a$]   (a) at (1.8,0) {};
  \node[circle,draw,inner sep=2pt,fill=gray!30,label=below:$v$] (v) at (3.6,0) {};
  \node[circle,draw,inner sep=2pt,fill=white,label=below:$b$]   (b) at (5.4,0) {};
  \node[circle,draw,inner sep=2pt,fill=gray!30,label=below:$w$] (w) at (7.2,0) {};
  \draw[very thick,red!70!black] (u) -- node[above]{$e_1$} (a);
  \draw (a) -- (v);
  \draw[very thick,red!70!black] (v) -- node[above]{$e_2$} (b);
  \draw (b) -- (w);
\end{tikzpicture}
\caption{A sharp example with $T=\{u,v,w\}$ and $F=\{e_1,e_2\}$. Here $\dist{G}{u}{v}=\dist{G}{v}{w}=2$, $\dist{G}{u}{w}=4$, and all three terminal-terminal distances become infinite in $G\setminus F$.}
\label{fig:semi-cut-example}
\end{figure}

The bound is sharp. For $k\ge 2$, let $G$ be the path on vertices $v_0,v_1,\dots,v_{2k-2}$, let $T=\{v_0,v_2,\dots,v_{2k-2}\}$, and let $F=\bigl\{\{v_{2i},v_{2i+1}\}: 0\le i\le k-2\bigr\}$. Then $|F|=|T|-1$, and every component of $G\setminus F$ contains exactly one terminal. Hence each terminal-terminal distance is finite in $G$ and infinite in $G\setminus F$. Figure~\ref{fig:semi-cut-example} shows the case $k=3$.

\begin{remark}\label{rem:trees}
If $G$ is a tree and $\emptyset\neq T\subseteq V$, then $T$ admits a semi-cutting set of size exactly $|T|-1$; by Theorem~\ref{thm:main} this is optimal. Indeed, fix $t_0\in T$ and root $G$ at $t_0$. For $t\in T\setminus\{t_0\}$, let $e_t$ be the first edge of the unique $t$--$t_0$ path, starting from $t$, and put $F\coloneqq\{e_t:t\in T\setminus\{t_0\}\}$. Each $e_t$ joins $t$ to its parent, so these edges are pairwise distinct and $|F|=|T|-1$. Given distinct $s,s'\in T$, let $x$ be the vertex closest to the root on the unique $s$--$s'$ path. At least one of $s,s'$ differs from $x$; after interchanging them if necessary, assume that $s\neq x$. Then $e_s$ lies on the $s$--$s'$ path. Hence $F$ meets every shortest path between distinct terminals, and Lemma~\ref{lem:geodesic-characterization} applies.
\end{remark}

Viewed through shortest paths, the problem asks for an edge transversal of the family of all shortest paths with endpoints in $T$. When $|T|=2$, this is a shortest-path disruption problem for one prescribed pair of terminals. When $|T|>2$, however, it is neither the classical multiterminal cut problem, which seeks total separation of the terminals, nor a vertex-covering problem for maximal geodesics in the whole graph. Rather, it occupies a distinct intermediate position: the family is defined metrically by terminal pairs, but the transversal is taken in edges. In Section~\ref{sec:related}, we discuss several closely related problems from the literature, highlighting their connections to and differences from the semi-cutting problem studied here.

The key reformulation is Lemma~\ref{lem:geodesic-characterization} below, which identifies semi-cutting with hitting all shortest paths whose endpoints lie in $T$ and serves as the pivot for both reductions. If some terminal is incident with an edge of $F$, then one removes that terminal from $T$ together with one incident edge from $F$; any shortest path between the remaining terminals that uses that edge still contains another edge of $F$. If no terminal is incident with $F$, then contracting the closed neighborhood of any terminal to a new terminal preserves the property while reducing the number of vertices. This dichotomy leads to a lexicographic induction on $(|T|,|V|)$.

\section{Preliminaries}

Let $G=(V,E)$ be a finite simple graph. A \emph{walk} from $u$ to $w$ is a sequence $P=(v_0,\dots,v_k)$ with $v_0=u$, $v_k=w$, and $\{v_{i-1},v_i\}\in E$ for $1\le i\le k$; its length is $\len{P}\coloneqq k$. A \emph{path} is a walk with pairwise distinct vertices. The distance $\dist{G}{u}{w}$ is the minimum length of a $u$--$w$ path when such a path exists, and $\infty$ otherwise. A path of length $\dist{G}{u}{w}$ is a \emph{shortest} path. For a path $P=(v_0,\dots,v_k)$, write $V(P)\coloneqq \{v_0,\dots,v_k\}$, $E(P)\coloneqq \bigl\{\{v_{i-1},v_i\}: 1\le i\le k\bigr\}$, and $P[v_i,v_j]\coloneqq (v_i,\dots,v_j)$ for $0\le i\le j\le k$. For $F\subseteq E$, let $G\setminus F\coloneqq (V,E\setminus F)$ and $V(F)\coloneqq \{v\in V:\ v\in e \text{ for some } e\in F\}$. For $v\in V$, let $N(v)\coloneqq \{w\in V:\ \{v,w\}\in E\}$ and $N[v]\coloneqq N(v)\cup\{v\}$.

We use the following vertex-set contraction operation. When the contracted vertex set induces a connected subgraph, this agrees with the usual contraction of a connected subgraph in minor theory.

\begin{definition}[Vertex-set contraction]\label{def:contraction}
Let $\emptyset\neq W\subseteq V$, and let $w^\ast\notin V$ be a new vertex. Define $\pi_W\colon V\to (V\setminus W)\cup\{w^\ast\}$ by $\pi_W(v)=w^\ast$ for $v\in W$ and $\pi_W(v)=v$ for $v\notin W$. The contracted graph $G/W$ has vertex set $(V\setminus W)\cup\{w^\ast\}$ and edge set
\[
E(G/W)=\bigl\{\{\pi_W(x),\pi_W(y)\}: \{x,y\}\in E,\ \pi_W(x)\neq \pi_W(y)\bigr\}.
\]
Thus $G/W$ is obtained by contracting $W$ to $w^\ast$; since we remain in the category of simple graphs, loops are deleted and parallel edges are identified. If $e=\{x,y\}\in E$ and $\pi_W(x)\neq \pi_W(y)$, we write $e/W\coloneqq \{\pi_W(x),\pi_W(y)\}\in E(G/W)$. For $F\subseteq E$, we also put $F/W\coloneqq \bigl\{\{\pi_W(x),\pi_W(y)\}: \{x,y\}\in F,\ \pi_W(x)\neq \pi_W(y)\bigr\}\subseteq E(G/W)$.
\end{definition}

The map
\[
\{e=\{x,y\}\in E : \pi_W(x)\neq \pi_W(y)\} \longrightarrow E(G/W),
\qquad e\longmapsto e/W,
\]
is well defined and surjective, although it need not be injective.

\begin{definition}[Semi-cutting set]\label{def:semi-cut}
Let $T\subseteq V$. An edge set $F\subseteq E$ is a \emph{semi-cutting set} for $T$ in $G$ if $\dist{G}{u}{w}<\dist{G\setminus F}{u}{w}$ for every two distinct vertices $u,w\in T$.
\end{definition}

\noindent We record two routine facts about walks and shortest paths that will be used repeatedly.

\begin{lemma}[Walks can be shortened to paths]\label{lem:walk-to-path}
Every walk from $u$ to $w$ contains a $u$--$w$ path of length at most the length of the walk.
\end{lemma}

\begin{proof}
Whenever the walk visits some vertex twice, deleting the closed subwalk between the two occurrences preserves the endpoints and does not increase the length. Iterating this operation yields a walk with no repeated vertices, hence a path.
\end{proof}

\begin{lemma}[Subpaths of shortest paths]\label{lem:subpath-shortest}
If $P=(v_0,\dots,v_k)$ is a shortest path in $G$ and $0\le i<j\le k$, then $P[v_i,v_j]$ is a shortest path from $v_i$ to $v_j$.
\end{lemma}

\begin{proof}
If $P[v_i,v_j]$ were not shortest, replacing it with a shorter $v_i$--$v_j$ path would produce a $v_0$--$v_k$ walk shorter than $P$. Lemma~\ref{lem:walk-to-path} would then give a shorter $v_0$--$v_k$ path, contrary to the choice of $P$.
\end{proof}

\begin{lemma}[Projection of walks under contraction]\label{lem:projected-walk}
Let $\emptyset\neq W\subseteq V$, and let $P=(v_0,\dots,v_k)$ be a walk in $G$. Apply $\pi_W$ to each vertex of $P$ and then delete every consecutive repetition. The resulting sequence is a walk in $G/W$ from $\pi_W(v_0)$ to $\pi_W(v_k)$, and its length is at most $k$.
\end{lemma}

\begin{proof}
Write $z_i=\pi_W(v_i)$ for $0\le i\le k$. Whenever $z_{i-1}\neq z_i$, Definition~\ref{def:contraction} gives $\{z_{i-1},z_i\}\in E(G/W)$. Deleting consecutive repetitions from $(z_0,\dots,z_k)$ therefore yields a walk in $G/W$ with the same endpoints. Only deletions are performed, so the length does not increase.
\end{proof}

\begin{lemma}[Geodesic characterization]\label{lem:geodesic-characterization}
Let $X\subseteq V$ and $M\subseteq E$. Then $M$ is a semi-cutting set for $X$ in $G$ if and only if, for every two distinct vertices $u,w\in X$, the distance $\dist{G}{u}{w}$ is finite and every shortest $u$--$w$ path in $G$ contains an edge of $M$.
\end{lemma}

\begin{proof}
Assume first that $M$ is semi-cutting for $X$. For distinct $u,w\in X$, the distance $\dist{G}{u}{w}$ cannot be infinite, because Definition~\ref{def:semi-cut} would then give $\infty<\dist{G\setminus M}{u}{w}$. If a shortest $u$--$w$ path in $G$ were disjoint from $M$, it would remain a $u$--$w$ path in $G\setminus M$, contradicting the strict inequality in Definition~\ref{def:semi-cut}.

Conversely, assume that every distinct $u,w\in X$ have finite distance in $G$ and that every shortest $u$--$w$ path in $G$ contains an edge of $M$. Fix such a pair $u,w$. If $u$ and $w$ are disconnected in $G\setminus M$, then $\dist{G\setminus M}{u}{w}=\infty>\dist{G}{u}{w}$. Otherwise a shortest $u$--$w$ path $Q$ in $G\setminus M$ is also a $u$--$w$ path in $G$, and equality $\len{Q}=\dist{G}{u}{w}$ is impossible because $Q$ is disjoint from $M$. Hence $\dist{G\setminus M}{u}{w}=\len{Q}>\dist{G}{u}{w}$.
\end{proof}

Lemma~\ref{lem:geodesic-characterization} will be used repeatedly as the geodesic reformulation of semi-cutting and is the common starting point of both reductions in the next section.

\section{Structural Reductions}

\begin{lemma}\label{lem:contraction-image-size}
Let $\emptyset\neq W\subseteq V$ and let $F\subseteq E$. Then $|F/W|\le |F|$.
\end{lemma}

\begin{proof}
The restriction of the map $e\mapsto e/W$ to $\{e=\{x,y\}\in F : \pi_W(x)\neq \pi_W(y)\}$ is surjective onto $F/W$. Hence $|F/W|\le |F|$.
\end{proof}

\begin{lemma}\label{lem:outside-edges}
Let $t\in V$, let $G^\prime\coloneqq G/N[t]$, and denote by $t^\prime$ the new vertex arising from the contraction. If $x,y\in V\setminus N[t]$ are distinct, then $\{x,y\}\in E(G^\prime)$ if and only if $\{x,y\}\in E(G)$; in that case $\{x,y\}/N[t]=\{x,y\}$.
\end{lemma}

\begin{proof}
If $\{x,y\}\in E(G)$, then Definition~\ref{def:contraction} gives $\{x,y\}\in E(G^\prime)$ and $\{x,y\}/N[t]=\{x,y\}$. Conversely, if $\{x,y\}\in E(G^\prime)$, then some edge $\{u,v\}\in E(G)$ satisfies $\{\pi_{N[t]}(u),\pi_{N[t]}(v)\}=\{x,y\}$. Since neither endpoint is $t^\prime$, both $u$ and $v$ lie outside $N[t]$, hence $\{u,v\}=\{x,y\}$.
\end{proof}

\begin{lemma}[Lifting shortest paths through a closed-neighborhood contraction]\label{lem:lifting-endpoint}
Let $t\in V$, let $G^\prime\coloneqq G/N[t]$, and denote by $t^\prime$ the new vertex arising from the contraction. Let $z\in V\setminus N[t]$, and let $P^\prime=(t^\prime,x_1,\dots,x_q=z)$ be a shortest path from $t^\prime$ to $z$ in $G^\prime$. Then there exists a vertex $v\in N(t)$ such that $P=(t,v,x_1,\dots,x_q=z)$ is a shortest path from $t$ to $z$ in $G$.
\end{lemma}

\begin{proof}
Because $P^\prime$ is a path, none of $x_1,\dots,x_q$ equals $t^\prime$; equivalently, all of them lie outside $N[t]$. The edge $\{t^\prime,x_1\}$ of $G^\prime$ therefore comes from an edge $\{v,x_1\}\in E(G)$ with $v\in N[t]$. Since $x_1\notin N[t]$, we must have $v\in N(t)$. Lemma~\ref{lem:outside-edges} shows that every edge $\{x_{i-1},x_i\}$ also belongs to $G$, so $P=(t,v,x_1,\dots,x_q)$ is a $t$--$z$ path in $G$ of length $\len{P^\prime}+1$.

Now let $R=(r_0=t,r_1,\dots,r_\ell=z)$ be any $t$--$z$ path in $G$. Since $z\notin N[t]$, we have $\ell\ge 2$ and $r_1\in N(t)$. Applying Lemma~\ref{lem:projected-walk} to the tail $(r_1,\dots,r_\ell)$ gives a $t^\prime$--$z$ walk in $G^\prime$ of length at most $\ell-1$, hence a $t^\prime$--$z$ path of length at most $\ell-1$ by Lemma~\ref{lem:walk-to-path}. Thus $\dist{G^\prime}{t^\prime}{z}\le \ell-1=\len{R}-1$. As $R$ was arbitrary, every $t$--$z$ path in $G$ has length at least $\dist{G^\prime}{t^\prime}{z}+1$, so $\dist{G}{t}{z}\ge \dist{G^\prime}{t^\prime}{z}+1$. Since $P$ is a $t$--$z$ path, we also have $\dist{G}{t}{z}\le \len{P}=\len{P^\prime}+1=\dist{G^\prime}{t^\prime}{z}+1$.
Hence $\len{P}=\dist{G}{t}{z}$, and $P$ is shortest.
\end{proof}

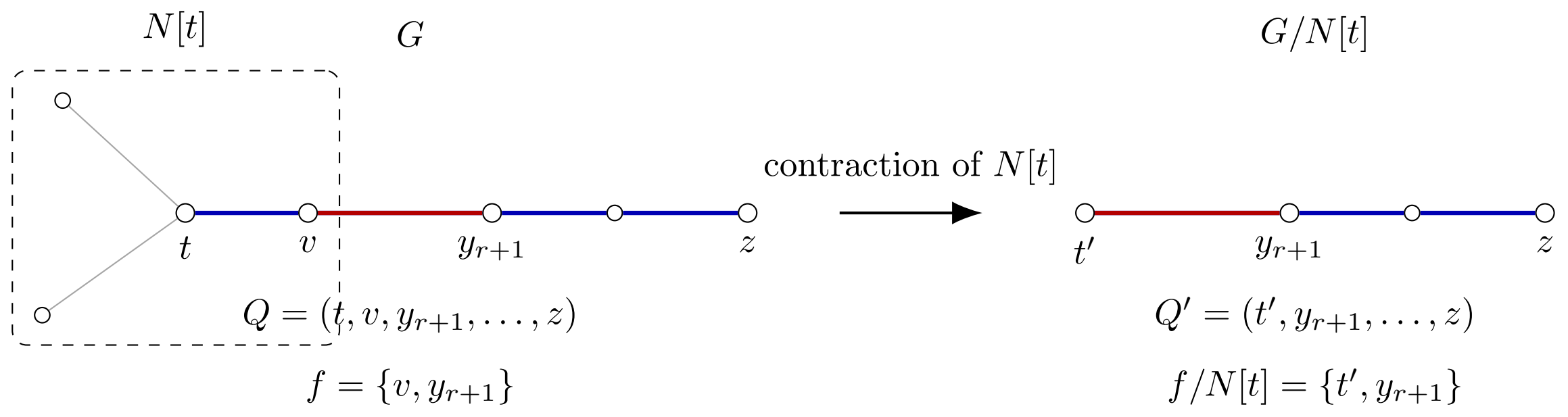
\begin{figure}[t]
\centering
\begin{tikzpicture}[scale=1, every node/.style={font=\small}]
  \begin{scope}[xshift=0cm]
    \node[circle,draw,inner sep=1.5pt,fill=white] (a) at (-0.2,1.1) {};
    \node[circle,draw,inner sep=1.5pt,fill=white] (b) at (-0.4,-1.0) {};
    \node[circle,draw,inner sep=1.8pt,fill=white,label=below:$t$] (t) at (1.0,0.0) {};
    \node[circle,draw,inner sep=1.8pt,fill=white,label=below:$v$] (v) at (2.2,0.0) {};
    \node[circle,draw,inner sep=1.8pt,fill=white,label=below:$y_{r+1}$] (y1) at (4.0,0.0) {};
    \node[circle,draw,inner sep=1.5pt,fill=white] (y2) at (5.2,0.0) {};
    \node[circle,draw,inner sep=1.8pt,fill=white,label=below:$z$] (z) at (6.5,0.0) {};
    \draw[gray!70] (a)--(t) (b)--(t) (t)--(v);
    \draw[very thick,blue!70!black] (t)--(v)--(y1)--(y2)--(z);
    \draw[very thick,red!70!black] (v)--(y1);
    \node[draw,dashed,rounded corners,fit=(a)(b)(t)(v),inner sep=6pt] {};
    \node at (0.9,1.8) {$N[t]$};
    \node at (3.2,-1.0) {$Q=(t,v,y_{r+1},\dots,z)$};
    \node at (3.2,-1.7) {$f=\{v,y_{r+1}\}$};
    \node at (3.2,1.75) {$G$};
  \end{scope}

  \draw[-{Latex[length=3mm]},thick] (7.4,0) -- (8.8,0);
  \node at (8.1,0.45) {\small contraction of $N[t]$};

  \begin{scope}[xshift=9.8cm]
    \node[circle,draw,inner sep=1.8pt,fill=white,label=below:$t^\prime$] (tp) at (0.0,0.0) {};
    \node[circle,draw,inner sep=1.8pt,fill=white,label=below:$y_{r+1}$] (yp1) at (2.0,0.0) {};
    \node[circle,draw,inner sep=1.5pt,fill=white] (yp2) at (3.2,0.0) {};
    \node[circle,draw,inner sep=1.8pt,fill=white,label=below:$z$] (zp) at (4.5,0.0) {};
    \draw[very thick,blue!70!black] (tp)--(yp1)--(yp2)--(zp);
    \draw[very thick,red!70!black] (tp)--(yp1);
    \node at (2.25,-1.0) {$Q^\prime=(t^\prime,y_{r+1},\dots,z)$};
    \node at (2.25,-1.7) {$f/N[t]=\{t^\prime,y_{r+1}\}$};
    \node at (2.25,1.75) {$G/N[t]$};
  \end{scope}
\end{tikzpicture}
\caption{Lifting a shortest $t^\prime$--$z$ path through the contraction in Proposition~\ref{prop:contraction}.}
\label{fig:contraction}
\end{figure}

\begin{proposition}[Contraction at a terminal not incident with $F$]\label{prop:contraction}
Let $T\subseteq V$ and $F\subseteq E$. Assume that $F$ is a semi-cutting set for $T$ in $G$. Let $t\in T\setminus V(F)$. Let $G^\prime\coloneqq G/N[t]$, and denote by $t^\prime$ the new vertex arising from the contraction. Put $T^\prime\coloneqq (T\setminus\{t\})\cup\{t^\prime\}$ and $F^\prime\coloneqq F/N[t]$. Then $F^\prime$ is a semi-cutting set for $T^\prime$ in $G^\prime$.
\end{proposition}

\begin{proof}
No terminal of $T\setminus\{t\}$ lies in $N[t]$. Indeed, if $u\in T\setminus\{t\}$ belonged to $N[t]$, then $\{t,u\}$ would be the unique shortest $t$--$u$ path in $G$, so Lemma~\ref{lem:geodesic-characterization} would force $\{t,u\}\in F$, contradicting $t\notin V(F)$. Thus $T\setminus\{t\}\subseteq V\setminus N[t]$. In particular, every terminal in $T\setminus\{t\}$ survives the contraction as the same vertex, no two terminals are identified, and $|T^\prime|=|T|$.

Every two distinct terminals of $T^\prime$ have finite distance in $G^\prime$. If $u,w\in T\setminus\{t\}$, take a $u$--$w$ path in $G$; if one terminal is $t^\prime$ and the other is $z\in T\setminus\{t\}$, take a $t$--$z$ path in $G$. In both cases, Lemmas~\ref{lem:projected-walk} and~\ref{lem:walk-to-path} produce a path in $G^\prime$ joining the corresponding terminals. By Lemma~\ref{lem:geodesic-characterization}, it remains to show that every shortest path in $G^\prime$ between distinct terminals of $T^\prime$ meets $F^\prime$.

Let $u,w\in T^\prime$ be distinct, and let $P^\prime$ be a shortest $u$--$w$ path in $G^\prime$.

\smallskip\noindent
\textbf{Case 1: $t^\prime\notin V(P^\prime)$.}
Then every vertex of $P^\prime$ lies outside $N[t]$, so $P^\prime$ is also a path in $G$ by Lemma~\ref{lem:outside-edges}. It is shortest in $G$. Indeed, if $R$ were a shorter $u$--$w$ path in $G$, then Lemma~\ref{lem:projected-walk} would project $R$ to a $u$--$w$ walk in $G^\prime$ of length at most $\len{R}<\len{P^\prime}$. Shortening that walk to a path by Lemma~\ref{lem:walk-to-path} would still give a $u$--$w$ path in $G^\prime$ of length strictly smaller than $\len{P^\prime}$, contradicting the choice of $P^\prime$. Since now $u,w\in T\setminus\{t\}\subseteq T$, Lemma~\ref{lem:geodesic-characterization} yields an edge $g\in E(P^\prime)\cap F$. Both ends of $g$ lie outside $N[t]$, hence $g/N[t]=g\in F^\prime$ by Lemma~\ref{lem:outside-edges}.

\smallskip\noindent
\textbf{Case 2: $t^\prime\in V(P^\prime)$.}
Because $P^\prime$ is a path, the vertex $t^\prime$ appears on it at most once. By symmetry, write $P^\prime=(y_0,\dots,y_m)$ with $y_r=t^\prime$ for some $r<m$. Put $z\coloneqq y_m\in T^\prime\setminus\{t^\prime\}=T\setminus\{t\}$ and $Q^\prime\coloneqq P^\prime[y_r,y_m]$. Since $Q^\prime$ is a path, all vertices $y_{r+1},\dots,y_m$ lie outside $N[t]$. By Lemma~\ref{lem:subpath-shortest}, $Q^\prime$ is a shortest $t^\prime$--$z$ path in $G^\prime$. Lemma~\ref{lem:lifting-endpoint} therefore provides a vertex $v\in N(t)$ such that $Q=(t,v,y_{r+1},\dots,y_m)$ is a shortest $t$--$z$ path in $G$; see Figure~\ref{fig:contraction}.

Since $t,z\in T$, Lemma~\ref{lem:geodesic-characterization} gives $E(Q)\cap F\neq\emptyset$. The first edge $\{t,v\}$ is not in $F$ because $t\notin V(F)$. Hence some edge $f\in E(Q)\cap F$ lies after the first step of $Q$. If $f=\{v,y_{r+1}\}$, then $\pi_{N[t]}(v)=t^\prime\neq y_{r+1}=\pi_{N[t]}(y_{r+1})$, so $f/N[t]=\{t^\prime,y_{r+1}\}$ is defined and is the first edge of $Q^\prime$. Otherwise $f=\{y_i,y_{i+1}\}$ for some $i\in\{r+1,\dots,m-1\}$, and both endpoints of $f$ lie outside $N[t]$, so Lemma~\ref{lem:outside-edges} gives $f/N[t]=f\in E(Q^\prime)$. Because $G^\prime$ is simple, distinct edges of $G$ may acquire the same image after contracting $N[t]$; however, the definition of $F^\prime=F/N[t]$ guarantees that any edge of $G^\prime$ which is the image of some edge of $F$ belongs to $F^\prime$, regardless of any additional preimages outside $F$. Thus the image of $f$ lies in $E(Q^\prime)\cap F^\prime\subseteq E(P^\prime)\cap F^\prime$.

Thus every shortest path in $G^\prime$ between distinct terminals of $T^\prime$ meets $F^\prime$. Lemma~\ref{lem:geodesic-characterization} completes the proof.
\end{proof}

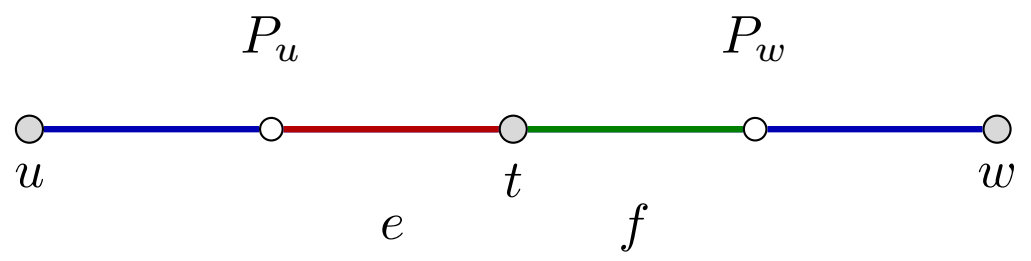
\begin{figure}[t]
\centering
\begin{tikzpicture}[scale=1, every node/.style={font=\small}]
  \node[circle,draw,inner sep=1.8pt,fill=gray!30,label=below:$u$] (u) at (0,0) {};
  \node[circle,draw,inner sep=1.5pt,fill=white] (x) at (1.6,0) {};
  \node[circle,draw,inner sep=1.8pt,fill=gray!30,label=below:$t$] (t) at (3.2,0) {};
  \node[circle,draw,inner sep=1.5pt,fill=white] (y) at (4.8,0) {};
  \node[circle,draw,inner sep=1.8pt,fill=gray!30,label=below:$w$] (w) at (6.4,0) {};
  \draw[very thick,blue!70!black] (u)--(x)--(t)--(y)--(w);
  \draw[very thick,red!70!black] (x)--(t);
  \draw[very thick,green!50!black] (t)--(y);
  \node at (1.6,0.6) {$P_u$};
  \node at (4.8,0.6) {$P_w$};
  \node at (2.4,-0.65) {$e$};
  \node at (4.0,-0.65) {$f$};
\end{tikzpicture}
\caption{In Lemma~\ref{lem:terminal-deletion}, if a shortest $u$--$w$ path uses the edge $e$ incident with $t$, then the opposite subpath still contains an edge of $F\setminus\{e\}$.}
\label{fig:split-path}
\end{figure}

\begin{lemma}[Deleting one terminal together with one incident edge]\label{lem:terminal-deletion}
Let $T\subseteq V$ and $F\subseteq E$. Assume that $F$ is a semi-cutting set for $T$ in $G$. Let $t\in T$, and let $e\in F$ be an edge incident with $t$. Put $T^-\coloneqq T\setminus\{t\}$ and $F^-\coloneqq F\setminus\{e\}$. Then $F^-$ is a semi-cutting set for $T^-$ in $G$.
\end{lemma}

\begin{proof}
If $|T^-|\le 1$, there is nothing to prove. Assume $|T^-|\ge 2$. By Lemma~\ref{lem:geodesic-characterization}, every two distinct vertices of $T^-$ have finite distance in $G$. It is therefore enough to show that every shortest path between two distinct vertices of $T^-$ contains an edge of $F^-$.

Take distinct $u,w\in T^-$, and let $P$ be a shortest $u$--$w$ path in $G$. Since $u,w\in T$, Lemma~\ref{lem:geodesic-characterization} yields $E(P)\cap F\neq\emptyset$. If $e\notin E(P)$, then already $E(P)\cap F^-\neq\emptyset$. Assume now that $e\in E(P)$. Write $P=(v_0,\dots,v_m)$ and let $v_j=t$. Because $u,w\neq t$, we have $0<j<m$. Put $P_u\coloneqq (v_0,\dots,v_j)$ and $P_w\coloneqq (v_j,\dots,v_m)$. By Lemma~\ref{lem:subpath-shortest}, these are shortest paths for the pairs $(u,t)$ and $(t,w)$, respectively.

Since $e$ is incident with $t$ and $t$ is an internal vertex of the path $P$, the edge $e$ lies on exactly one of $P_u$ and $P_w$; thus the choice of the other subpath is unambiguous. Let $R$ denote that other subpath. Then $R$ is a shortest path with endpoints $t$ and one of $u,w$, and $e\notin E(R)$. As both endpoints of $R$ lie in $T$, Lemma~\ref{lem:geodesic-characterization} gives an edge $f\in E(R)\cap F$. Since $e\notin E(R)$, we have $f\neq e$, so $f\in F^-$. Because $E(R)\subseteq E(P)$, the path $P$ contains an edge of $F^-$. Therefore every shortest path between distinct terminals of $T^-$ meets $F^-$, and Lemma~\ref{lem:geodesic-characterization} yields the claim.
\end{proof}

\section{A Sharp Lower Bound}

\begin{theorem}\label{thm:main}
Let $G=(V,E)$ be a finite simple undirected graph, let $T\subseteq V$, and let $F\subseteq E$. If $F$ is a semi-cutting set for $T$ in $G$, then $|F|\ge |T|-1$.
\end{theorem}

\begin{proof}
For $(k,n)\in \mathbb{N}_0^2$, let $\mathcal{S}(k,n)$ denote the statement that every simple graph $H$ on $n$ vertices, every terminal set $X\subseteq V(H)$ with $|X|=k$, and every semi-cutting set $M\subseteq E(H)$ for $X$ satisfy $|M|\ge k-1$. We prove $\mathcal{S}(k,n)$ for all $(k,n)$ by induction on $(k,n)$ in the lexicographic order, where $(a,b)\prec(c,d)$ means that either $a<c$, or $a=c$ and $b<d$.

There is nothing to prove when $k\le 1$. Fix $(k,n)$ with $k\ge 2$, assume that $\mathcal{S}(k',n')$ holds for every $(k',n')\prec (k,n)$, and let $G=(V,E)$ be a simple graph on $n$ vertices with a terminal set $T$ of size $k$ and a semi-cutting set $F\subseteq E$. By Lemma~\ref{lem:geodesic-characterization}, every terminal pair in $T$ has finite distance in $G$, so both reductions below are metrically meaningful even without any assumption that $G$ be connected.

We distinguish two cases according to whether $T\cap V(F)$ is empty.

\smallskip\noindent
\textbf{Case 1: $T\cap V(F)\neq\emptyset$.}
Choose $t\in T\cap V(F)$. Then some edge $e\in F$ is incident with $t$. Put $T^-\coloneqq T\setminus\{t\}$ and $F^-\coloneqq F\setminus\{e\}$. Lemma~\ref{lem:terminal-deletion} shows that $F^-$ is semi-cutting for $T^-$ in $G$. As $(|T^-|,|V(G)|)\prec (k,n)$, induction gives $|F^-|\ge k-2$, and therefore $|F|=|F^-|+1\ge k-1$.

\smallskip\noindent
\textbf{Case 2: $T\cap V(F)=\emptyset$.}
Choose $t\in T$. Since $k\ge 2$, pick $u\in T\setminus\{t\}$. Lemma~\ref{lem:geodesic-characterization} gives $\dist{G}{t}{u}<\infty$, so $N(t)\neq\emptyset$ and therefore $|N[t]|\ge 2$. Let $G^\prime\coloneqq G/N[t]$, denote by $t^\prime$ the new vertex arising from the contraction, and set $T^\prime\coloneqq (T\setminus\{t\})\cup\{t^\prime\}$ and $F^\prime\coloneqq F/N[t]$. By Proposition~\ref{prop:contraction}, the set $F^\prime$ is semi-cutting for $T^\prime$ in $G^\prime$. Moreover, $|T^\prime|=k$ and $|V(G^\prime)|=|V|-|N[t]|+1\le n-1$. Hence $(|T^\prime|,|V(G^\prime)|)\prec (k,n)$, so $|F^\prime|\ge k-1$ by induction. Lemma~\ref{lem:contraction-image-size} gives $|F^\prime|\le |F|$, and therefore $|F|\ge k-1$.

Together, the terminal-deletion step in Case~1 and the contraction step in Case~2 close the induction and prove the lower bound.
\end{proof}

\section{Related Problems and Further Directions}\label{sec:related}

When $|T|=2$, Lemma~\ref{lem:geodesic-characterization} turns semi-cutting into the problem of hitting all shortest paths between one prescribed pair of terminals. This places the parameter near single-pair shortest-path disruption and most-vital-edge or most-vital-arc questions. Corley and Sha~\cite{corley-sha}, as well as Malik, Mittal, and Gupta~\cite{malik-mittal-gupta}, study the deletion of a prescribed number of edges, vertices, or arcs so as to maximize the increase of the shortest distance between one fixed source-target pair in a weighted network, while Nardelli, Proietti, and Widmayer~\cite{nardelli-proietti-widmayer} analyze the most vital edge of a shortest source-target path together with the corresponding replacement shortest paths after removing one edge of that path. Our parameter is different even in this two-terminal regime: it does not optimize a distance increase, but asks whether a given edge set meets \emph{all} shortest paths between the prescribed terminals; for larger $T$, the relevant family comprises the shortest paths for every terminal pair.

\begin{figure}[t]
\centering
\begin{tikzpicture}[scale=1, every node/.style={font=\small}]
  \node[circle,draw,inner sep=1.8pt,fill=gray!30,label=left:$u$] (u) at (0,0) {};
  \node[circle,draw,inner sep=1.8pt,fill=gray!30,label=right:$v$] (v) at (4,0) {};
  \node[circle,draw,inner sep=1.8pt,fill=gray!30,label=below:$w$] (w) at (1.1,-3) {};
  \node[circle,draw,inner sep=1.5pt,fill=white,label=above:$a$] (a) at (1.4,1.2) {};
  \node[circle,draw,inner sep=1.5pt,fill=white,label=right:$b$] (b) at (3.1,-1.6) {};
  \node[circle,draw,inner sep=1.5pt,fill=white,label=left:$c$] (c) at (-0.3,-1.5) {};
  \node[circle,draw,inner sep=1.5pt,fill=white,label=below:$x$] (x) at (1.5,-1.0) {};

  \draw (u)--(a);
  \draw[very thick,red!70!black] (a)-- node[above right]{$e_1$} (v);
  \draw (v)--(b);
  \draw[very thick,red!70!black] (b)-- node[right]{$e_2$} (w);
  \draw (w)--(c);
  \draw[very thick,red!70!black] (c)-- node[left]{$e_3$} (u);
  \draw (a)--(x);
  \draw (b)--(x);
  \draw (c)--(x);
\end{tikzpicture}
\caption{A semi-cutting set that is not a multiterminal cut. Here $T=\{u,v,w\}$ and $F=\{e_1,e_2,e_3\}$. In $G$, every terminal pair is at distance $2$; in $G\setminus F$, every terminal pair is still connected but now at distance $4$.}
\label{fig:not-multiterminal}
\end{figure}
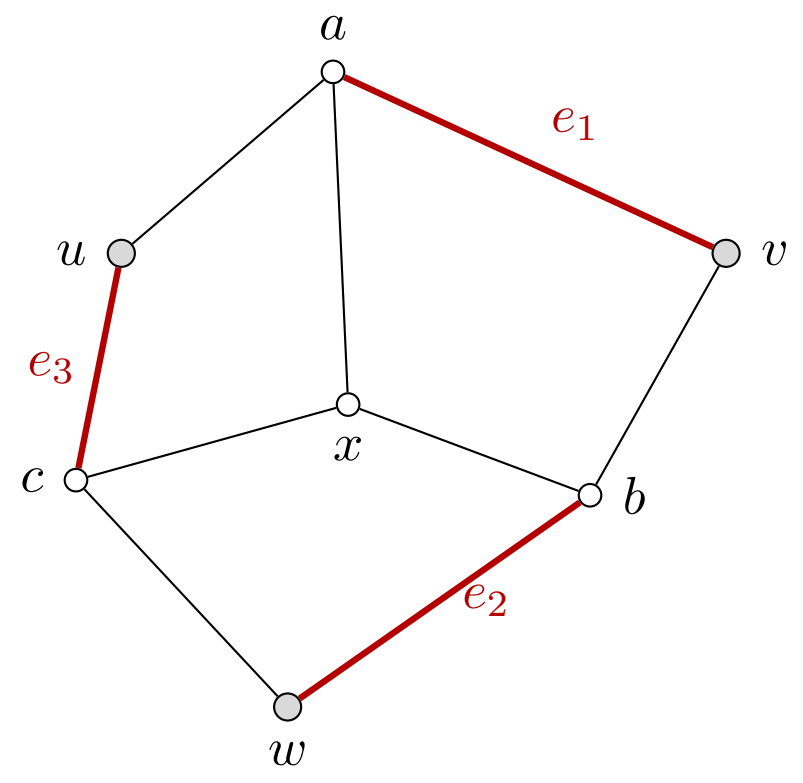

For larger terminal sets, the closest classical separation problem is the multiterminal cut problem of Dahlhaus, Johnson, Papadimitriou, Seymour, and Yannakakis~\cite{dahlhaus1994}. There the aim is to delete edges so that each terminal is separated from all the others, and isolating cuts and multiterminal cuts are defined precisely in that sense. Semi-cutting imposes a strictly weaker topological requirement but a subtler metric one: every terminal-terminal distance must increase, yet the terminals may remain in the same connected component after deleting $F$. Figure~\ref{fig:not-multiterminal} illustrates exactly this phenomenon, showing that semi-cutting is not reducible to multiterminal separation.

The closest conceptual relatives of the semi-cutting problem are vertex parameters defined by geodesic families. Harary, Loukakis, and Tsouros~\cite{harary-loukakis-tsouros} introduce geodetic closure and the geodetic number, while Oellermann and Peters-Fransen~\cite{oellermann-peters-fransen} use intervals $I[u,v]$ in their study of strong metric dimension; both belong to the broader tradition of shortest-path-based vertex parameters. Closer to the present setting, Peterin and Semani\v{s}in~\cite{peterin-semanisin} study a minimum vertex set meeting all maximal shortest paths. Manuel, Bre\v{s}ar, and Klav\v{z}ar~\cite{manuel-bresar-klavzar-transversal} define the geodesic-transversal number as a minimum vertex set hitting all maximal geodesics, and in a subsequent paper~\cite{manuel-bresar-klavzar-packing} they introduce the dual geodesic packing parameter. Semi-cutting differs from this line in three respects: it is edge-based rather than vertex-based; it concerns the terminal-defined family rather than maximal geodesics of the whole graph; and it requires hitting all shortest terminal-terminal paths rather than maximal geodesics.

Two problems arise naturally from Theorem~\ref{thm:main}. One is to characterize the equality case $|F|=|T|-1$. The other is to study minimum semi-cutting sets as edge transversals of the family of shortest paths with endpoints in $T$.

Theorem~\ref{thm:main} thus isolates a purely structural phenomenon, independent of algorithmic considerations: once semi-cutting is viewed through the terminal-defined family of shortest terminal-terminal paths, contraction and terminal deletion force a universal lower bound. Further progress therefore seems likely to come from structural characterizations of extremal configurations rather than from procedures.

\section*{Acknowledgments}

We thank Marcin Anholcer, Grzegorz Gutowski, Micha{\l} Laso{\'n}, T{\'a}ssio Naia, Jakub Przyby{\l}o, Micha{\l} Tuczy{\'n}ski, Alexander Wolff, Mariusz Zaj{\k{a}}c, and Anna Zych-Pawlewicz for fruitful discussions, and Davide Bil{\`o} for bringing this problem to our attention.

\end{document}